\documentclass[11pt]{amsart}
\usepackage{amssymb,amsthm,amsmath}
\newcommand{\M}{\mathcal{M}}

\newcommand{\RR}{\mathbb{R}}

\newcommand{\A}{\mathcal{A}}

\newcommand{\cD}{\mathcal{D}}

\newcommand{\E}{\mathcal{E}}
\newcommand{\D}{\partial}

\newcommand{\<}{\langle}

\newcommand{\loc}{\mathrm{loc}}
\newcommand{\intr}{\mathrm{int}}
\newcommand{\vol}{\mathrm{vol}}
\renewcommand{\div}{\mathrm{div}}
\renewcommand{\>}{\rangle}
\newtheorem{thm}{Theorem}
\newtheorem{lem}[thm]{Lemma}

\theoremstyle{definition}

\title{The Calder\'on problem is an inverse source problem}
\author{Jan Cristina.}
\thanks{The author was supported by the Finnish Cultural Foundation with the grant ``Harmonic maps, coordinate gauges, and anisotropic inverse problems''}
\date{}
\begin{document}
\maketitle
\begin{abstract}
  We prove that  uniqueness for the Calder\'on problem on a Riemannian manifold with boundary  follows from a hypothetical unique continuation property for the elliptic operator $\Delta+V+(\Lambda^{1}_{t}-q)\otimes (\Lambda^{2}_{t}-q)$ defined on $\D\M^{2}\times[0,1]$ where $V$ and $q$ are potentials and $\Lambda^{i}_{t}$ is a Dirichlet--Neumann operator at depth $t$. This is done by showing that the difference of two Dirichlet--Neumann maps is equal to the Neumann boundary values of the solution to an inhomogeneous equation for said operator, where the source term is a measure supported on the diagonal of $\D\M^{2}$.
\end{abstract}
\section{Introduction}
Calder\'on's inverse boundary value problem asks whether the Cauchy data at the boundary of an elliptic second order pseudo-differential operator $\div(\sigma d)$ determine the coefficients $\sigma$.  It has been solved in a great deal of generality on Euclidean domains under the assumption of scalar coefficients $\sigma:\Omega\to \RR^{+}$ \cite{Calderon_in_plane,Sylvester_Uhlmann,Lipschitz_conductivities,CaroRogers_Lipschitz_conductivities}.  

The problem becomes significantly more difficult when the coefficients $\sigma$ are assumed to be anisotropic, that is given by a symmetric positive-definite tensor-field.  Nonetheless there is a very strong suggestion of general uniqueness provided by the proof in the case of real-analytic conductivities \cite{Lassas_Uhlmann_Taylor}.  The concept of limiting Carleman weights was impressively used to prove uniqueness of the Calder\'on problem for metrics in a conformal class, provided the class admitted such a weight \cite{Limiting_Carleman_weights}, but the condition is rather limiting on the geometry and topology of the spaces under consideration \cite{obstructions_to_carleman_weights},\cite{Non_Carleman_weight_manifolds}.  

In the present paper I show that Calder\'on's problem on a closed manifold with boundary $\M$ can be reduced to studying a related unique continuation problem on the space $\D\M\times\D\M\times[0,\varepsilon]$.

I introduce the following unique continuation property:  A second order differential operator $\A:C^{\infty}(\D\M^{2}\times[0,1])\to C^{\infty}(\D\M^{2}\times[0,1])$ is said to have the \emph{off diagonal unique continuation property} (ODUCP) if for any $u\in \cD'(\D\M^{2}\times[0,1])$ $\A u=0$ on $(\D\M^{2}\setminus \Delta)\times[0,1]$, $u|_{0}=0$ and $\D_{t} u|_{0}=0$ implies that $u$ equals $0$ on $\D\M^{2}\times [0,\varepsilon]$. 

Given a $C^{5}$-smooth compact  manifold with boundary $\M$, endowed with a $C^{4}$ smooth Riemannian metric $g$, we can consider the Dirichlet--Neumann map $\Lambda_{g}$ which takes a function $f:\D\M\to \RR$ and maps it to $\D_{\nu}u$, where $u$ is the solution to the boundary value problem
\begin{equation}
  \Delta_{g}u=0,\quad u|_{\D\M}=f.
  \label{eqn:Laplace_bvp}
\end{equation}
For $n\geq 3$ consider the associated problem for a scalar multiple of $g$, $\gamma g$ where $\gamma:\M\to \RR$ is a smooth function bounded away from zero.  This can be identified with the boundary value problem for a Schr\"odinger equation \cite{Sylvester_Uhlmann}. If $\Delta_{\gamma g}u=0$, then setting $\sigma=\gamma^{n/2-1}$ the function $\sigma^{1/2}u$ solves the boundary value problem
\[(\Delta_{g} +Q)\sigma^{1/2}u=0 \quad \sigma^{1/2}u|_{\D\M}=\sigma^{1/2}f.\]
where $Q=\sigma^{-1/2}\Delta_{g}\sigma^{1/2}$.  
Consequently the Dirichlet--Neumann map associated to $\gamma g$ can be equated with that for the Schr\"odinger equation
\[\Lambda_{\gamma g}=\sigma^{1/2}\Lambda_{Q}\sigma^{-1/2} +\D_{\nu}\sigma^{1/2}.\]
Consequently given two weights $\gamma_{1}$ and $\gamma_{2}$, they define the same Dirichlet--Neumann maps if and only if the Dirichlet--Neumann maps for the associated Schr\"odinger equations are equal, and the weights are equal at the boundary.  Let $Q^{1}$ and $Q^{2}$ denote the associated potentials, and let $\Lambda^{1}$ and $\Lambda^{2}$ denote the associated Dirichlet--Neumann maps.

Given a Riemannian manifold with boundary and $C^{2}$ metric tensor, we can define Fermi or boundary normal coordinates via the map $\Psi:\D\M\times[0,\varepsilon]\to \M$ $(x,t)\mapsto \gamma_{x,-\nu}(t)$, where $\gamma_{x,-\nu}$ is the inward normal oriented geodesic starting at $x$ at time $t$.  For $\varepsilon$ sufficiently small this is a diffeomorphism onto its image \cite{Inverse_boundary_spectral_problems}.  

In these coordinates, the metric $g$ takes the special form 
\[g=dt^{2}+h_{t},\]
where $h_{t}$ is a metric tensor on $\D\M$.  We let $d\vol_{h_{t}}=e^{\mu_{t}}d\vol_{h_{0}}$.

Then we can define a family of operators $\Lambda^{i}_{t}$ for $t\in[0,\varepsilon]$ to be the associated Dirichlet--Neumann operators for the  operators $\Delta_{g}+Q^{i}$ restricted to the submanifold $\M\setminus \psi(\D\M\times[0,t))$.  Now we can define a metric on $\D\M\times\D\M\times[0,\varepsilon]=$ by
\[G=dt^{2}+h_{t,x}+h_{t,y}.\]  The associated volume form is given by $dt\wedge d\vol_{t}=dt\wedge d\vol_{0}e^{\mu_{t}(x)+\mu_{t}(y)}$ for $(x,y)\in \D\M\times\D\M$.
  
With this formalism we are able to introduce our associated operator $\A:C^{\infty}(\D\M^{2}\times[0,\varepsilon]\to C^{\infty}(\D\M^{2}\times[0,\varepsilon])$ for the Calder\'on problem:
  \begin{equation}
    \A=\Delta_{G}-\dot{\mu}_{t}(x)\dot{\mu}_{t}(y)+Q^{1}(t,x)+Q^{2}(t,y)+ (\Lambda^{1}_{t}-\dot{\mu}_{t}(x))\otimes (\Lambda^{2}_{t}-\dot{\mu}_{t}(y))
    \label{eqn:A_operator_definition}
  \end{equation}
  We call $\A$ the evolution squared operator, because it arises as the product of an evolution operator and it's adjoint (\emph{cf} \S \ref{sec:evo_squared}).  We are thus able to condition uniqueness for the Calder\'on problem on the ODUCP of $\A$:
\begin{thm}
  \label{thm:uniq_cont_loc_cald}Let $\A$ be the operator defined in \eqref{eqn:A_operator_definition}.  If $\A$ has the ODUCP, then $\Lambda^{1}=\Lambda^{2}$ if and only if $Q^{1}=Q^{2}$ on $\D\M\times[0,\varepsilon]$ and $\Lambda^{1}_{\varepsilon}=\Lambda^{2}_{\varepsilon}$.
\end{thm}
The proof of Theorem \ref{thm:uniq_cont_loc_cald} follows by constructing a solution to a boundary value problem for $\A$ on $\D\M^{2}\times[0,\varepsilon]$ 
\begin{align}
\nonumber  &\A \varphi= \Phi\\
  &\varphi(0)=0\quad\text{and}\quad [\D_{t}+\Lambda^{1}_{\varepsilon}\otimes I+I\otimes \Lambda^{2}_{\varepsilon}]\varphi(\varepsilon)=\Lambda^{2}_{\varepsilon}-\Lambda^{1}_{\varepsilon}
  \label{eqn:A_bvp}
\end{align}
where $\Phi$ is the measure supported on $\Delta\times [0,1]$ given by 
\[\Phi(f)=\int_{\D\M\times[0,1]}f(x,x,t)[Q^{1}(t,x)-Q^{2}(t,x)]\:d\vol_{h_{t}}(x)\:dt.\]
This special solution will have $t$-derivative at $t=0$ equal to $\Lambda^{1}-\Lambda^{2}$ as a distribution on $\D\M\times\D\M$.  Consequently it follows from the ODUCP that if $\Lambda^{1}-\Lambda^{2}=0$ then $\varphi=0$.  However, elementary arguments can be used to show that if $Q^{1}-Q^{2}\neq 0$ then $\varphi\neq 0$ on $(\D\M^{2}\setminus \Delta)\times [0,\varepsilon]$.

Subsequently using the existence of an exhaustion (Lemma \ref{lem:exists_exhaustion}) this can be turned into the following global contrapositive of the preceding theorem.
\begin{thm}
  \label{thm:uniq_cont_impl_cald_uniq} 
  Suppose $\Lambda^{1}=\Lambda^{2}$, then $Q^{1}-Q^{2}\neq 0$  only if there is a collar neighbourhood $U$ of $\D\M$ such that the operator $\A$ associated to some Fermi coordinates of the complement of $U$ does not have the ODUCP.
\end{thm}
\section{Preliminaries}
  We work on a $C^{k+1}$-smooth  compact manifold $\M$ with boundary, and assume it has a $C^{k}$-smooth metric $g$.  The necessity of such a high degree of smoothness arises from Fermi coordinates.

  Let $\Psi:\D\M\times[0,\varepsilon)\to \M$ be the mapping taking $x,t$ to the inward normal oriented geodesic starting at $x$ at time $t$.  $\Psi$ is $C^{k-1}$ diffeomorphism for a $C^{k}$ metric. Consequently the pullback metric $\Psi^{*}g$ on $\D\M\times[0,\varepsilon]$ is $C^{k-2}$. Under the diffeomorphism the metric takes the form
  \[dt^{2}+h_{t}\]
  where $h_{t}$ is a metric tensor on $\D\M$ \cite{Inverse_boundary_spectral_problems}. Let $\Sigma_{t}$ denote $\M\setminus \Psi(\D\M\times[0,t))$.

We define Sobolev spaces on our manifolds via a smooth (that is as smooth as the manifold allows) partition of unity $\psi_{i}$, supported on a set $U_{i}$ with a coordinate chart $\varphi_{i}:U_{i}\to \RR^{n}_{+}$.  A measurable function $u$ is in $H^{s}(\M)$ if $\psi_{i}u\circ \varphi_{i}$ is in $H^{s}(\RR^{n}_{+})$ for every $s$.  If $\M$ is $C^{k}$ then this is valid for $s\in [0,k]$. 

For $s\in [0,k-1]$ we can define the spaces $H^{s}(\M,\Lambda^{l}\M)$ to be the space of measurable $l$-forms $\alpha$ for which $\varphi_{i}^{*}\psi_{i}\alpha\in H^{s}(\RR^{n}_{+},\Lambda^{l}\RR^{n})$.

We define the spaces $H^{s}(\D\M)$ and $H^{s}(\D\M,\Lambda^{l}\D\M)$ similarly.  We define the space $H^{-s}(\D\M,\Lambda^{l}\D\M)$ to be the dual of $H^{s}(\D\M, \Lambda^{n-1-l}\D\M)$ for $s\in[0,k-1]$.  This negates the need for a volume form if it is undesirable.  If $l=0$ then this definition can be extended to $s\in[0,k]$.  A norm $\|\cdot\|_{s}$ is fixed, although it is not particularly important which one.  For instant $\<(I+\Delta)^{s}\cdot,\cdot\>$ where $\Delta$ is some fixed Laplace--Beltrami operator.

  Given a potential $Q:\M\to \RR$, consider the equation
  \begin{align}
    &\Delta u+Qu=0\nonumber\\
    &u|_{\D\Sigma_{t}}\circ \Psi(t,\cdot)=f,
    \label{eqn:harmonic_function}
  \end{align}
  where $\Delta$ denotes the Laplace--Beltrami operator for  $(\M,g)$.
  We define the map $\E^{t}_{Q}:H^{1/2}(\D\M)\to H^{1}(\Sigma_{t})$ to be the \emph{solution operator} for equation \eqref{eqn:harmonic_function}, so $\E_{Q}(f)=u$.
  The Dirichlet--Neumann map $\Lambda_{Q}^{t}:H^{1/2}(\D\M)\to H^{-1/2}(\D\M)$ takes $f$ to the the normal derivative $\D_{\nu}\E_{Q}^{t}(f)\circ \Psi(t,\cdot)$.

This can be formulated elegantly as 
\[\int_{\D\M}\Lambda^{t}_{Q}(f_{1})(f_{2})\:d\vol_{t}=\int_{\Sigma_{t}}d\E_{Q}^{t}(f)\wedge\star_{g}d\E_{Q}^{t}(f_{2})+\int_{\Sigma_{t}}\E_{Q}^{t}(f_{1})Q\E_{Q}^{t}(f_{2})\:d\vol,\]
where $\star_{g}$ is the Hodge-star operator associated to $g$, $d\vol_{0}$ is the Volume form associated to $h_{0}$ on $\D\M$, and $d\vol$ is the volume form associated to $g$. Of course an important point is that
\[\int_{\D\M}\Lambda^{t}_{Q}(f_{1})(f_{2})\:d\vol_{t}=\int_{\Sigma_{t}}d\E_{Q}^{t}(f)\wedge\star_{g}d\tilde{\E}^{t}(f_{2})+\int_{\Sigma_{t}}\E_{Q}^{t}(f_{1})Q\tilde{\E}^{t}(f_{2})\:d\vol,\]
where $\tilde{\E}$ is an arbitrary extension operator, because the difference $\tilde{\E}^{t}-\E_{Q}^{t}\in H^{1}_{0}(\Sigma_{t})$.

Given two potentials $Q^{1}$ and $Q^{2}$ on $\M$, we can consider the difference in the Dirichlet--Neumann maps, here denoted for parsimony's sake by $\Lambda_{1}$ and $\Lambda_{2}$ respectively.  By a standard integration by parts technique we can show that
\[\int_{\D\M}(\Lambda_{1}-\Lambda_{2})(f_{1})f_{2}\:d\vol_{0}=\int_{\M}\E_{1}(f_{1})(Q^{1}-Q^{2})\E_{2}(f_{2})\:d\vol.\]
By splitting the right hand integral, we can express this as
\begin{multline}
  \label{eqn:layer_strip_decomposition}
  \int_{\D\M}(\Lambda_{1}-\Lambda_{2})(f_{1})f_{2}\:d\vol_{0}=\\
  \int_{0}^{\varepsilon}\int_{\D\M}\E_{1}(f_{1})(x,t)(Q^{1}-Q^{2})(x,t)\E_{2}(f_{2})(x,t)\:d\vol_{t}\:dt\\
  \quad+\int_{\D\M}(\Lambda^{\varepsilon}_{1}-\Lambda^{\varepsilon}_{2})(\E_{1}(f_{1})(\varepsilon,\cdot)\E_{2}(f_{2})(\varepsilon,\cdot)\:d\vol_{\varepsilon}.
\end{multline}
\section{The tautological evolution equation}
A key observation for the results herein is that the map $t\mapsto (x\mapsto \E_{Q}(f)(t,x))$ is the solution to an evolution equation
\begin{equation}
  \D_{t}u_{t}=-\Lambda^{t}_{Q}(u_{t}). \qquad u_{0}=f
  \label{eqn:tautalogical_evolution_equation}
\end{equation}
This is the \emph{tautological evolution equation} for the boundary value problem \eqref{eqn:harmonic_function}.  This observation follows from the fact that $\D_{t}=-\D_{\nu}$ for the manifold $\Sigma_{t}$ for all $t\in [0,\varepsilon]$ in Fermi coordinates.

Given two different potentials, and functions $u$ and $v$ satisfying
\[\D_{t}u_{t}=-\Lambda^{1}_{t}u_{t}\qquad\text{and}\qquad \D_{t}v_{t}=-\Lambda^{2}_{t}v_{t}\]
respectively we will be testing the difference in potentials agains their product
\[\int_{\D\M\times[0,\varepsilon)}u_{t}(x)(Q^{1}-Q^{2})(t,x)v_{t}(x)\:d\vol_{t}\:dt.\]
Now the product $u_{t}(x) v_{t}(x)$ does not \emph{a priori} satisfy an evolution equation like \eqref{eqn:tautalogical_evolution_equation}, but the tensor product $u_{t}\otimes v_{t}$ define as $(x,y)\in\D\M^{2}\mapsto u_{t}(x)v_{t}(y)$ satisfies the evolution equation
\begin{equation}
  \D_{t}(u_{t}\otimes v_{t})=-A_{t}u_{t}\otimes v_{t},
  \label{eqn:tensor_evolution}
\end{equation}
where  $A_{t}=\Lambda^{1}_{t}\otimes I+I\otimes \Lambda^{2}_{t}$.
\begin{lem}Let $u_{t}\otimes v_{t}$ satisfy
  \[\D_{t}(u_{t}\otimes v_{t})=-A_{t}(u_{t}\otimes v_{t}),\]
  then
  \label{lem:evolved_tensor_prod_solves_A}
  \[[\Delta_{G}+Q^{1}(t,x)+Q^{2}(t,y)-\dot{\mu}(x)\dot{\mu}(y)/2+2(\Lambda^{1}_{t}-\dot{\mu}/2)\otimes (\Lambda^{2}_{t}-\dot{\mu}/2)]u_{t}\otimes v_{t}=0\]
\end{lem}
\proof  The operator $\Delta_{G}= \Delta_{g}^{x}+\Delta_{g}^{y}-\D_{t}^{2}-[\dot{\mu}(x)+\dot{\mu}(y)]\D_{t}$.  When we apply this to $u_{t}\otimes v_{t}$ we arrive at
\begin{align*}
  \Delta_{G}u_{t}\otimes v_{t}&=\Delta_{g}^{x}u_{t}\otimes v_{t}+u_{t}\otimes \Delta_{g}^{y}v_{t}-\D_{t}^{2}u_{t}\otimes v_{t}-u_{t}\otimes \D_{t}^{2}v_{t}\\
  &\quad -2\D_{t}u_{t}\otimes \D_{t}v_{t}-\dot{\mu}(x)\D_{t}u_{t}\otimes v_{t}-u_{t}\otimes \dot{\mu}(y)\D_{t} v_{t} \\
  &\quad -\dot{\mu}(x) u_{t}\otimes \D_{t} v_{t}-\dot{\mu}(y)\D_{t}u_{t}\otimes v_{t}\\
  &=-Q^{1}(t,x)u_{t}\otimes v_{t}-u_{t}\otimes Q^{2}(t,y) v_{t}-2\Lambda^{1}_{t}(u_{t})\otimes\Lambda^{2}_{t}(v_{t})\\
 &\quad +\dot{\mu}(x)u_{t}\otimes\Lambda^{2}(v_{t})+\dot{\mu}(y)\Lambda^{1}(u_{t})\otimes v_{t},\\
\end{align*}
where we have replaces $\D_{t}u_{t}$ with $-\Lambda_{t}^{1}u_{t}$ and likewise for $\D_{t}v_{t}$.
\begin{lem}
  \label{lem:DN_map_satisfies_ext_props}
  Suppose $g$ and $Q$ are $C^{k}$-smooth, $k\geq 3$, then $t\mapsto \Lambda^{t}_{Q}$ is weakly $C^{1}$ i.e. for every $u$ and $v\in H^{1/2}(\D\M)$, the map
  \[t\mapsto \Lambda^{t}_{Q}(u)(v)\]
  is $C^{1}$.  Furthermore
  \[\D_{t}\Lambda^{t}_{Q}-(\Lambda^{t}_{Q})^{2}=-\Delta_{h_{t}}-Q\circ\Psi(t,\cdot) -\dot{\mu}_{t}\Lambda^{t}_{Q},\]
and
\begin{equation}
  \label{eqn:DN_lower_bound}
\<(1-\Delta)^{s}\Lambda^{t}_{Q} u,u\>_{t}\geq C_{1}\|u\|_{s+1/2}^{2}-C_{2}\|u\|_{s}^{2}
\end{equation}
for $s\in[-k,k]$
\end{lem}
\proof
The proof of this follows from layer stripping arguments in \cite{Sylvester_layer_stripping}, which were also applied in \cite{riem_inc_EIT}.
\begin{align*}
  \int_{\D\M}(\Lambda_{Q}^{t+h}(u)(v)&\:d\vol_{t+h}-\Lambda^{t}_{Q}(u)(v))\:d\vol_{t}\\
  &=\int_{\Sigma_{t+h}}d\E^{t+h}_{Q}(u)\wedge\star_{g}d\E^{t+h}_{Q}(v)-\int_{\Sigma_{t}}d\E_{Q}^{t}(u)\wedge\star_{g}d\E_{Q}^{t}(v)\\
  &\quad+\int_{\Sigma_{t+h}}\E^{t+h}_{Q}(u)Q\E^{t+h}_{Q}(v)\:d\vol-\int_{\Sigma_{t}}\E^{t}_{Q}(u)Q\E^{t}_{Q}(v)\:d\vol\\
  &=\int_{\Sigma_{t+h}}\bigg[d\E_{Q}^{t+h}(u-\E_{Q}^{t}(u)|_{\Sigma_{t+h}})\wedge\star_{g}d\E_{Q}^{t+h}(v)\\
  &\qquad\qquad\qquad+d\E_{Q}^{t}(u)\wedge\star_{g}d\E_{Q}^{t+h}(v-\E_{Q}^{t}(v)|_{\Sigma_{t+h}})\bigg]\\
  &\quad+\int_{\Sigma_{t+h}}\bigg[\E_{Q}^{t+h}(u-\E_{Q}^{t}(u)|_{\Sigma_{t+h}})Q\E_{Q}^{t+h}(v)\\
  &\qquad\qquad\qquad+\E_{Q}^{t}(u)Q\E_{Q}^{t+h}(v-\E_{Q}^{t}(v)|_{\Sigma_{t+h}})\bigg]\:d\vol\\
  &\quad-\int_{\Sigma_{t}\setminus \Sigma_{t+h}}d\E_{Q}^{t}(u)\wedge\star_{g}d\E_{Q}^{t}(v)i-\int_{\Sigma_{t}\setminus\Sigma_{t+h}}\E^{t}_{Q}(u)Q\E^{t}_{Q}(v)\:d\vol.
\end{align*}
If we divide by $h$ and let $h$ tend to $0$. we arrive at 
\begin{align*}
  \D_{t}(\Lambda^{t}_{Q}\:d\vol_{t})=((\Lambda^{t}_{Q})^{2}-\Delta_{t}-Q)\:d\vol_{t},
\end{align*}
but $\D_{t}d\vol_{t}=\dot{\mu}_{t} d\vol_{t}$.
Given that the principle symbol of $(\Lambda^{t}_{Q})^{2}$ is equal to $\Delta_{t}$ \cite{riem_inc_EIT} \cite{Uhlmann_Lee}, this yields that the $\D_{t}\Lambda^{t}$ is a bounded operator $H^{s+1}\to H^{s}$ uniformly in $t$. For the lower bound see \cite{riem_inc_EIT}.
\endproof

\section{Evolution squared and the singular source}
\label{sec:evo_squared}
The second order equation also factorises as the product of two evolution equations in different directions which motivates the name \emph{evolution squared}.  This is in turn justified when we try to make sense of the distributional boundary value problem, which is greatly facilitated compared to heavier machinery, such as the Boutet de Monvel calculus \cite{Boutet_de_monvel_calculus}.   

\begin{lem} The evolution squared operator $\A$ is equal to 
  \[\A=(\D_{t}^{*}+A_{t})(\D^{t}+A_{t}).\]
\end{lem}
\proof First we note that $\D_{t}^{*}=-\D_{t}-\dot{\mu}_{t}(y)-\dot{\mu}_{t}(x)$.
Here we make use of the fact that 
\[\D_{t}\Lambda^{i}_{t}-(\Lambda^{i}_{t})^{2}=-\Delta_{h}-Q^{i}-\dot{\mu}\Lambda^{i}_{t},\]
from which we deduce that
\[A^{2}-\D_{t}A=\Delta_{h_{t}\times h_{t}}+Q^{1}(t,x)+Q^{2}(t,y)+\dot{\mu}(x)\Lambda^{1}_{t}\otimes I+I\otimes \dot{\mu}(y)\Lambda^{2}_{t}+2\Lambda^{1}_{t}\otimes \Lambda^{2}_{t}.\]
Of course
\[\Delta_{G}=\D_{t}^{*}\D_{t}+\Delta_{h_{t}\times h_{t}},\]
and $[\D_{t},A_{t}]=\D_{t}(A_{t})$.  When we expand
\[(\D_{t}^{*}+A_{t})(\D_{t}+A_{t})=\D_{t}^{*}\D_{t}+A_{t}^{2}-[\D_{t},A_{t}]-(\dot{\mu}(x)+\dot{\mu}(y))A_{t}\]
with these, we arrive at the desired result.
\endproof
Now  we can see how $\A$ defines a reasonable elliptic operator.
%
We define a distributional solution $\varphi\in \cD'(\D\M^{2}\times[0,\varepsilon],\D\M^{2}\times\{0,\varepsilon\})$ of \eqref{eqn:A_bvp} as one for which 
    \[
      \int_{\D\M^{2}\times[0,\varepsilon]}\varphi \A U\:d\vol=\int_{\D\M\times[0,\varepsilon]}(Q^{1}-Q^{2})(t,x)U(t,x,x)\:d\vol
    -(\Lambda^{1}_{\varepsilon}-\Lambda^{2}_{\varepsilon})(U(\varepsilon,\cdot,\cdot))\]
    for every smooth function $U$ equal to $0$ on $\D\M^{2}\times\{0\}$ and for which $\D_{t}U=-A_{\varepsilon} U$ for $t=\varepsilon$.  The space $\cD'(\D\M^{2}\times[0,\varepsilon], \D\M^{2}\times\{0,\varepsilon\})$ is defined to be the topological dual of $C^{\infty}(\D\M\times[0,\varepsilon])$.

\begin{lem}
\label{lem:evolution_squared_bvp_existence}
Assume $g$ and $Q^{i}$ are  in $C^{\lceil n/2\rceil}$ where $n=\dim \D\M$.  There is a distributional solution to 
\[\A\varphi=[Q^{1}(t,x)-Q^{2}(t,y)]\delta_{x,y}\]
\[\varphi(0,x,y)=0,\quad (\D_{t}+A_{\varepsilon})\varphi(\varepsilon,x,y)=(\Lambda^{1}_{\varepsilon}-\Lambda^{2}_{\varepsilon})\]
    \end{lem}
    \proof We start by solving the inhomogeneous evolution equation 
\[(\D_{t}^{*}+A_{t})\psi_{t}=[-\D_{t}-(\dot{\mu}(x)+\dot{\mu}(y))+A_{t}]\psi_{t}=0,\] 
subject to $\psi_{\varepsilon}=-(\Lambda^{1}_{\varepsilon}-\Lambda^{2}_{\varepsilon})$ as a distribution on $\D\M\times\D\M$.  It is an element of $H^{-1}(\D\M\times \D\M)$

We consider a solution $\psi\in H^{-1}$, equipped, with an inner product of the form 
\[\<f,f\>_{-1}=\<(I+\Delta)^{-1}(u),u\>,\]
where $\<\cdot,\cdot\>$ is an $L^{2}$ inner product on $\D\M\times \D\M$.  Then $A_{t}$ is monotone with respect to this inner product by virtue of $\eqref{eqn:DN_lower_bound}$.

Then we have an evolution triple $H^{-1}$, $H^{-1/2}$ $H^{-3/2}$ with this inner product, and so we have unique existence of a solution \cite[Theorem 23.A]{nonlinear_fnl_analysis_and_apps_IIa} in the space $H^{1}([0,1],H^{-1}(\D\M),H^{-1/2}(\D\M),H^{-3/2}(\D\M))$.

We then solve the inhomogeneous initial value problem
\begin{equation}
  [\D_{t}+\A_{t}]\hat{\varphi}=\psi(t), \quad\hat{\varphi}(0)=0.
\end{equation}

We also solve the inhomogeneous initial value problems
\[[-\D_{t}+\dot{\mu}(x)+\dot{\mu}(y)+A_{t}]\tilde{\psi}=\Phi_{t},\quad \tilde{\psi}(\varepsilon)=0,\]
and
\[[\D_{t}+A_{t}]\tilde{\varphi}=\tilde{\psi},\quad \tilde{\varphi}(0)=0.\]
However $\Phi_{t}$ is no longer in $H^{-1}$, but $H^{-s}$ for $s>n/2$ for every $t$.  So we introduce a new evolution triple $(H^{-s},\<(1+\Delta)^{-n/2}\cdot,\cdot\>),H^{1/2-s}$, $H^{-1/2-s}$.
Lastly we set $\varphi=\tilde{\varphi}+\hat{\varphi}$.
\endproof

We apply integration by parts to get the following:
\begin{lem}
  \label{lem:boundaryvalue}
  Let $\varphi$ be a distributional solution to \eqref{eqn:A_bvp}: 
  \[\A \varphi=\Phi,\quad \varphi(0)=0,\quad (\D_{t}+A_{t})\varphi(\varepsilon)=(\Lambda^{1}_{\varepsilon}-\Lambda^{2}_{\varepsilon}),\]
  then
  \[(\Lambda^{1}-\Lambda^{2})(u)(v)=\int_{\D\M^{2}}\D_{t}\varphi(0)u\otimes \:d\vol_{0}.\]
\end{lem}
\proof
Let $\psi=(\D_{t}+\A_{t})\varphi$. We begin with \eqref{eqn:layer_strip_decomposition} to get
\begin{align*}
  (\Lambda^{1}-\Lambda^{2})(u_{0})(v_{0})&=\int_{0}^{\varepsilon}\Phi(u_{t}\otimes v_{t})\:dt+(\Lambda^{1}_{\varepsilon}-\Lambda^{2}_{\varepsilon})(u_{\varepsilon}\otimes v_{\varepsilon})\\
  &=\int_{\D\M^{2}\times[0,\varepsilon]}(\D_{t}^{*}+A_{t})(\D_{t}+A_{t})\varphi \:u_{t}\otimes v_{t}\:d(\mu\otimes \mu)\:dt\\
  &\quad+(\Lambda^{1}_{\varepsilon}-\Lambda^{2}_{\varepsilon})(u_{\varepsilon}\otimes v_{\varepsilon})\\
  &=-\int_{\D\M^{2}}\psi(0)\:u_{0}\otimes v_{0}+\int_{\D\M^{2}}\psi(\varepsilon)\:u_{\varepsilon}\otimes v_{\varepsilon}\:d(\mu\otimes \mu)\\
  &\quad+(\Lambda^{1}_{\varepsilon}-\Lambda^{2}_{\varepsilon})(u_{\varepsilon}\otimes v_{\varepsilon})\\
  &=-\int_{\D\M^{2}}\psi(0)\:u_{0}\otimes v_{0}\:d(\mu\otimes\mu)
\end{align*}
Lastly we note that $\psi(0)=(\D_{t}+A_{0})\varphi(0)=\D_{t}\varphi(0)$ because $\varphi(0)=0$.
\endproof

We now need to show that $Q^{1}-Q^{2}\neq 0$ implies $\varphi\neq 0$ off the diagonal:
\begin{lem}
  \label{lem:phi_neq0_off_diag}
  Let $\Phi$ be the measure on $\D\M^{2}\times[0,\varepsilon]$ given by $\delta_{x,y}[Q^{1}_{t}(x)-Q^{2}_{t}(y)]$.  If $[Q^{1}-Q^{2}](x_{0},t_{0})\neq 0$ and 
  \[\A \varphi=\Phi\]
  then $\varphi\neq 0$ on  $(\D\M^{2}\setminus \Delta)\times[0,\varepsilon]$
\end{lem}
\proof
We will endeavor to show that $\varphi\in W_{\loc}^{1,p}$ where $p<\frac{n}{n-1}$, and $\varphi\not\in W_{\loc}^{1,p}$ for $p>\frac{n}{n-1}$, i.e. $\|d\varphi\|_{p}=\infty$.  First we note that $\Phi\in W^{-1,p}(\M)$ for $p<n/(n-1)$ and $\A$ is an elliptic second order pseudodifferential operator, so $\varphi\in W^{1,p}_{\loc}$.  To see that $\varphi\not\in W^{1,p}_{\loc}$ for $p=n/(n-1)$, consider a sequence of test functions $\varphi_{k}$ such that $\varphi(x,x,t)\to \infty$  for $x\in B(x_{0},r)$ and $t\in (t-\delta,t+\delta)$,  while $\|\varphi\|_{1,n}\leq C<\infty$. First we start with a function $\eta$ supported on $B(x_{0},r_{0})\times(t_{0}-\delta,t_{0}+\delta)$ such that
\[\int_{t_{0}-\delta}^{t_{0}+\delta}\int_{B(x_{0},r_{0}}\eta (x,t)(Q^{1}-Q^{2})(x,t)\:d\vol_{t}(x)\:dt\neq 0.\]  
Let  $\zeta_{k}=\min\{k,\log(\log(d_{t}(x,y)))\}\eta(x,t)$. 
Then
\begin{align*}
  \Phi(\zeta_{k})&=\int_{\D\M^{2}\times[0,\varepsilon]}(\D_{t}+A_{t})(\varphi)(D_{t}+A_{t})\xi_{k}\:d\vol_{t}\:dt\\
  &\leq \|\varphi\|_{1,n/(n-1)}\|\xi_{k}\|_{1,n}
\end{align*}
But $\Phi(\zeta_{k})\to \infty$ as $k\to \infty$, while $\|\zeta_{k}\|_{1,n}\leq C<\infty$.
\endproof

Armed with most of the necessary tools, we can now prove Theorem \ref{thm:uniq_cont_loc_cald}.
\proof[Proof ot Theorem \ref{thm:uniq_cont_loc_cald}]
By Lemma \ref{lem:boundaryvalue}, $(\Lambda^{1}-\Lambda^{2})=\D_{t}\varphi(0)$ where $\varphi$ is the solution to \eqref{eqn:A_bvp}.  By the hypothesized ODUCP if $\D_{t}\varphi(0)=0$ and $\varphi(0)=0$ then $\varphi=0$ everywhere.  But by Lemma \ref{lem:phi_neq0_off_diag} $\varphi$ cannot be zero if $Q^{1}-Q^{2}$ is non-zero.
\endproof

\section{Exhaustion and the proof of Theorem \ref{thm:uniq_cont_impl_cald_uniq}}

In order to prove Theorem \ref{thm:uniq_cont_impl_cald_uniq} we need a lemma for guaranteeing the existence of an exhaustion for smooth manifolds with boundary:
\begin{lem}
  \label{lem:exists_exhaustion}
  Let $\M$ be a $C^{2}$ smooth manifold with boundary.  There exists a map $\Phi:\D\M\times[0,1]\to \M$ which is a diffeomorphism onto its image when restricted to $\D\M\times[0,1)$, such that $\M\setminus\Phi(\D\M\times[0,1])$ is meager.  
\end{lem}
\proof  The proof requires the existence of a smooth triangulation $\Sigma$ of $\M$ \cite{smooth_triangulation}.  Let $\tau_{i}:\sigma_{n}\to \M$ be a smooth $n$-simplex, i.e. $C^{2}$ up to the boundary of each sub-simplex, a diffeomorphism on the interior.  By $\sigma_{n}$ we denote the set
\[\left\{x\in \RR^{n+1}:x^{i}\geq 0,\,\, \sum_{i=1}^{n+1}x^{i}=1\right\}.\]
We say $\sigma_{k}\subset \sigma_{l}$ by the natural inclusion of $\RR^{k+1}\subset \RR^{l+1}$ along the first $k+1$ components.

We will assume by induction that we have a diffeomorphism from some $N$ simplices $\Psi_{N}:\D\M\times[0,1]\to \bigcup_{i=1}^{N}\tau_{i}(\sigma_{n})$.  If there is a simplex  whose interior is disjoint from the image of $\Psi_{N}$, then there is a simplex whose interior is disjoint from the image of $\varphi$ and which neighbours some $\tau_{i}$.   Let $\tau_{N+1}$ denote this simplex, let $\nu$ denote a mutual facet of $\tau_{i}$ and $\tau_{N+1}$.  The goal will be to smoothly push $\Psi_{N}$ from $\tau_{i}$ through $\nu$ to $\tau_{N+1}$.  Let $\sigma_{n-1}$ denote the pre-image of $\nu$ in $\tau_{i}$.  Define the map $\zeta:\sigma_{n-1}\times[0,1)\to \sigma_{n}$ by 
  \[(x,t)\mapsto (x(1-t\eta(x)),t\eta(x)),\]
  where $\eta(x)=\prod_{i=1}^{n} x_{i}$.  Let $U=\zeta(\intr(\sigma_{n-1}))\times[0,1))$, and let $\sigma_{n}'$ denote the set $R_{n+1}(\sigma_{n})$ where $R_{n+1}:\RR^{n+1}\to \RR^{n+1}$ is the reflection $(x,x_{n+1})\mapsto (x,-x_{n+1})$.  $\sigma_{n}'\cap \sigma_{n}$ can be canonically smoothed.  We will construct a map $U\to U\cup \sigma_{n}'$: consider 
    \[ x\mapsto (x,0)+\min\{x_{i}\}(-\mathbf{1}/n,1)\quad x\in \sigma_{n-1}
      \]
      where $(\mathbf{1}=(1,1,\ldots,1)\in \RR^{n}$.  This takes $\sigma_{n-1}$ to the complement of $\sigma_{n-1}$ in $\D\sigma_{n}$.      We define $\min_{\varepsilon}(x_{1},\ldots, x_{n})$ inductively by
      \[\min_{\varepsilon}(x_{1},\ldots, x_{n})=\frac{1}{n}\sum_{i=1}\min_{\varepsilon}(\min_{\varepsilon}(x_{1},\ldots \hat{x}_{i},\ldots, x_{n}),x_{i}),\]
      where 
      \[\min_{\varepsilon}(x,y)=x+y+\sqrt{\varepsilon^{2}+|x-y|^{2}}-\varepsilon.\]
      Define
      \[\tilde{\zeta}(x,t)\mapsto \begin{cases}
	\zeta(x,\psi(x,t))& t\in[1/2,1]\\
	(x,0)+(1-2t)\min_{t}\{x_{i}\}(-\mathbf{1}/n),1)&t\in[0,1/2],
    \end{cases}\]
    where $\psi(t)$ maps $[1/2,1]$ to $[0,1]$ monotonically, and is the identity in a neighbourhood of $1$.     Lastly we must smooth $\tilde{\zeta}$ in a neighbourhood of $\sigma_{n-1}$. 

    Finally we define 
    \[\Psi_{N+1}(x,t)=\begin{cases}
	\Psi_{N}(x,t)& t<1\quad\Psi_{N}(x,t)\not\in \tau_{i}(U)\\
	\tau_{N+1}(\tilde{\zeta}(\zeta^{-1}(\tau_{i}^{-1}(\Psi_{N}(x,t)))))&\Psi_{N}(x,t)\in \tau_{i}(\zeta(\tilde(\zeta)^{-1}(\sigma_{n}'))))\\
	\tau_{i}(\tilde{\zeta}(\zeta^{-1}(\tau_{i}^{-1}\Psi_{N}(x,t))))&\Psi_{N}(x,t)\in \tau_{i}(U).
    \end{cases}\]

    To complete the induction we need an initial step.  Of course there is no diffeomorphism from $\D\M\to \tau_{1}$ for a single simplex, rather we must start with a collar.  Because $\M$ is $C^{2}$ we can apply the collar neighbourhood theorem, to yield a diffeomorphism from $\Psi_{0}:\D\M\times[0,1]\to \M$ onto its image.  Let $\Sigma=\M\setminus \Psi_{0}(\D\M\times[0,1)]$.  Choose a triangulation of $\Sigma$, and extend it to a triangulation of $\M$ by restricting the triangulation of $\Sigma$, to $\D\Sigma$, identifying $\D\M$ and $\D\Sigma$,  and canonically triangulating $\D\M\times[0,1]$ and joining this triangulation to that of $\Sigma$ to get a triangulation of $\M$.  Now $\Psi_{0}$ is  a diffeomorphism $\D\M\times[0,1]$ to some sub-triangulation of $\M$. 
\endproof
\proof[Proof of Theorem \ref{thm:uniq_cont_impl_cald_uniq}] 
We make use of Lemma \ref{lem:exists_exhaustion} to give us an exhaustion of $\M$ via the boundary $\Psi:\D\M\times[0,1]\to \M$.  If $(Q^{1}-Q^{2})(\Psi(x,t))=0$ for $t\in[0,s]$  then by \eqref{eqn:layer_strip_decomposition} $\Lambda^{1}_{s}-\Lambda^{2}_{s}=0$. Consequently by Theorem \ref{thm:uniq_cont_loc_cald} we can extend it two a boundary normal neighbourhood.  But then, because $\Psi$ is a homeomorphism, it also true for $t\in [0,s+\varepsilon]$ for some $\varepsilon$.  But then it is true for every $t\in [0,1)$ as the set on which it is true is closed and open.  Finally because $\M\setminus \Psi(\D\M\times[0,1])$ is meager, $Q^{1}=Q^{2}$ on all of $\M$.\endproof

\section{Some concluding remarks}
Although the off-diagonal unique continuation property is a strong assumption for the operator $\A$ there is some evidence to suggest it might hold.  Generically, unique continuation properties for pseudodifferential operators are not known and probably false, and it seems unlikely that an appropriate Carleman estimate could be derived for $\A$ because of the non-local behaviour of the Dirichlet--Neumann maps contained therein. However the work of Caffarelli and Silvestre \cite{Caffarelli_Silvestre} shows us that Dirichlet--Neumann maps have very strong unique continuation properties.  Our operator $\A$ is the sum of a differential operator and the tensor product of two Dirichlet--Neumann maps, so perhaps clever application of arguments like those in \cite{Caffarelli_Silvestre} could be used to derive such a unique continuation.  

Nonetheless, there should be no confusion that the ODUCP is a stronger condition than the uniqueness of the Calder\'on problem, however, if we restrict ourselves to the study of an operator $\A$ defined for $\Lambda^{1}_{t}=\Lambda^{2}_{t}$, then the ODUCP is equivalent to the uniqueness for the linearised Calder\'on problem for $\D\M\times[0,\varepsilon]$ with the metric $dt^{2}+h_{t}$.  Consequently it seems counterintuitive that uniqueness for the Calder\'on problem would be true, while the ODUCP would be false. 
\bibliographystyle{plain}
\def\cprime{$'$}

\end{document}